\input amstex 
\documentstyle{amsppt} 
%
%
\def\em#1{{\it #1\/}}
\def\tit#1{{\sl #1\/}}
%
%

\def\Integers{{\bold Z}}
\def\integers{\Integers}

\def\union{\cup}

\def\intersect{\cap}

\def\inverse{{-1}}

\def\suchthat{\bigm|}

\def\cross{\times}

\def\directsum{{\bold +}}
 
%
%
%
\def\func#1{\operatorname{#1}}

\def\interior{\func{int}}

\def\blocus{\raise0pt\hbox{$\Upsilon$}}

\def\CiteSection#1{\S #1}
\def\section#1{\CiteSection#1}

\def\nth#1{#1^{\hbox{\eightpoint th}}}


\def\cal{\Cal}
 
\def\qed{\hfill$\diamond$}
 




\magnification 1200
\hsize = 6.5 true in
 
\topmatter
 
\title\nofrills Central quotients of biautomatic groups \endtitle
 
\author Lee Mosher \endauthor
 
\address {Department of Mathematics and Computer Science,  Rutgers
University,  Newark, NJ 07102} \endaddress
 
\email mosher\@andromeda.rutgers.edu \endemail
 
\date 
\vbox{
\centerline{March 22, 1994}
}
\enddate
 
\abstract The quotient of a biautomatic group by a subgroup of the center is
shown to be biautomatic. The main tool used is the Neumann-Shapiro triangulation
of $S^{n-1}$, associated to a biautomatic structure on $\integers^n$. As an
application, direct factors of biautomatic groups are shown to be biautomatic.
\endabstract
 
\thanks 
The author was partially supported by NSF grant \# DMS-9204331 \endgraf  
This preprint is available from the Magnus archive on the American Mathematical
Society's e-MATH host computer. Send the message ``HELP'' to
mail-server\@e-math.ams.org
\endthanks
 
\endtopmatter
 
\document
 
\def\subgroup{<} 
\def\Star{\func{Star}} 

\def\Min{\func{Min}}
 
Biautomatic groups form a wide class of finitely presented groups with interesting
geometric and computational properties. These groups include all word hyperbolic
groups, all fundamental groups of finite volume Euclidean and hyperbolic
orbifolds, and all braid groups \cite{E$\ldots$}. A biautomatic group satisfies a
quadratic isoperimetric inequality, has a word problem solvable in quadratic
time, and has a solvable conjugacy problem. The class of biautomatic groups has
several interesting closure properties. For instance, the centralizer of a finite
subset of a biautomatic group is biautomatic, as is the center of the whole group
\cite{GS}. Also, biautomatic groups are closed under direct products
\cite{E$\ldots$}. The theory of biautomatic groups is briefly reviewed below. 
 
We present a technique for putting biautomatic structures on central quotients
of biautomatic groups:
 
\proclaim{Theorem A} Let $G$ be a biautomatic group, and let $C$ be a subgroup
of $\cal ZG$, the center. Then $G/C$ is biautomatic. \endproclaim
 
This result has several applications. Our first application answers a question
posed by Gersten and Short \cite{GS, cf. proposition 4.7}:
 
\proclaim{Theorem B} Direct factors of biautomatic groups are biautomatic.
\endproclaim
 
\demo{Proof} Suppose $G \cross H$ is biautomatic. The centralizer of $H$ is $G
\cross \cal ZH$, and this is a biautomatic group by \cite{GS, corollary 4.4}.
Then $\cal ZH$ is a subgroup of the center of $G \cross \cal ZH$, so by theorem
A, $G \cross \cal ZH / \cal ZH = G$ is biautomatic. \qed\enddemo
 
Several recent discoveries have pointed to the useful concept of \em{poison
subgroups}. For instance, the group $\integers^2$ is poison to word hyperbolic
groups: if a group contains a $\integers^2$ subgroup it cannot be word
hyperbolic. A wider class of subgroups poison to word hyperbolicity are those
which have an infinite index central $\integers$ subgroup \cite{CDP, corollaire
7.2}. Our next theorem says that for biautomatic groups, this class of poison
subgroups completely collapses to $\integers^2$:
 
\proclaim{Theorem C} If the biautomatic group $G$ contains a subgroup with an
infinite index central $\integers$ subgroup, then $G$ contains a $\integers^2$
subgroup. \endproclaim 
 
\demo{Proof} The hypothesis says that $G$ has an infinite cyclic subgroup $Z$ of
infinite index in its centralizer $\cal C_Z$. The group $\cal C_Z$ is biautomatic
by \cite{GS, corollary 4.4}. Since $Z$ is central in $\cal C_Z$, then by theorem
A the group $\cal C_Z / Z$ is biautomatic. This group is infinite, so by
\cite{E$\ldots$, example 2.5.12} it has an element of infinite order. Any
infinite cyclic subgroup of $\cal C_Z / Z$ pulls back to a $\integers^2$
subgroup of $\cal C_Z \subgroup G$. \qed\enddemo
 
Gersten and Short ask whether a biautomatic group can have an infinitely
generated abelian subgroup \cite{GS, p. 154}. We can reduce this problem as
follows:
 
\proclaim{Theorem D} Suppose there is a biautomatic group with an infinitely
generated abelian subgroup. Then either there is a biautomatic group with an
infinite rank abelian subgroup, or there is a biautomatic group which is not
virtually torsion free. \endproclaim
 
\demo{Proof} Suppose the biautomatic group $G$ has an abelian subgroup $H$,
infinitely generated and of finite rank $n \ge 0$. If $n \ge 1$, choose an element
$h \in H$ of infinite order. By \cite{GS, corollary 4.4}, the centralizer $\cal
C_h$ of $h$ in $G$ is biautomatic. Let $\cal Z \cal C_h$ be the center of $\cal
C_h$. By theorem A, $\cal C_h/\cal Z \cal C_h$ is biautomatic. Note that $h \in
\cal Z \cal C_h \subgroup H \subgroup \cal C_h$, so the image of $H$ in $\cal
C_h/\cal Z \cal C_h$ is an infinitely generated, abelian subgroup of rank $\le
n-1$. By induction, we obtain a biautomatic group $\Gamma$ with an infinitely
generated abelian subgroup $A$ of rank 0, i.e. an infinite abelian torsion
subgroup. If $\Gamma$ had a torsion free subgroup $K$ of finite index, then
there would be two elements $a \ne b \in A$ such that $1 \ne b^\inverse a \in K$,
hence $b^\inverse a$ has infinite order; but $b^\inverse a \in A$ has finite
order. \qed\enddemo
 
\subhead Proof of theorem A \endsubhead
 
First we reduce theorem A to a special case:
 
\proclaim{Theorem E} If $G$ is a biautomatic group and $Z \subgroup G$ is an
infinite cyclic central subgroup, then $G/Z$ is biautomatic. \endproclaim
 
\demo{Proof of theorem A} Let $G$ be a biautomatic group, and let $C$ be a
subgroup of the center $\cal ZG$. Since $\cal ZG$ is biautomatic it is finitely
generated. Hence $C$ is a finite rank central subgroup, say of rank $k \ge 0$.
Now peel off factors of $\integers$ one at a time, as in the proof of \cite{GS,
proposition 4.7}. If $k \ge 1$ let $Z$ be any infinite cyclic subgroup of $C$.
Applying theorem E it follows that $G/Z$ is biautomatic, and $C/Z$ is a central
subgroup of rank $k-1$. Repeating this argument $k$ times, we see that there is a
finite index free abelian subgroup $C' \subgroup C$ such that $C/C'$ is a central
finite subgroup in the biautomatic group $G/C'$. But the quotient of any
biautomatic group by any finite normal subgroup is easily seen to be biautomatic;
projecting the biautomatic structure from the total group to the quotient group
gives a biautomatic structure on the quotient. Hence, $G/C = (G/C')/(C/C')$ is
biautomatic. \qed\enddemo
 
The remainder of the paper is devoted to proving theorem E. 
 
\subhead Review of biautomatic groups \endsubhead
 
An \em{alphabet} is a finite set. A \em{word} over an alphabet $\cal A$ is a
finite sequence of elements in $\cal A$. The empty word is sometimes denoted
$\epsilon$. The set of all words in $\cal A$ is denoted $\cal A^*$, and this
forms a monoid under the operation of concatenation, with $\epsilon$ as the
identity. A \em{language over $A$} is a subset of $\cal A^*$.
 
The length of a word $w$ is denoted $\ell_w$. If $w$ is written in the
form $w = w_1 w_2$ then $w_1$ is called a \em{prefix} subword and $w_2$ is a
\em{suffix} subword. If $w$ is written as $w = w_1 w_2 w_3$ then $w_2$ is called
an \em{infix} subword, with associated prefix $w_1$ and suffix $w_3$. For any
integer $t \ge 0$, $w(t)$ denotes the prefix subword of $w$ of length $t$ if $t
\le \ell_w$, and $w(t)=w$ otherwise. 
 
We adopt the graph theoretic notion of a \em{finite state automaton} over an
alphabet $\cal A$. This is a finite directed graph $M$ whose vertices are called
\em{states}, together with a labelling of each edge by a letter of $\cal A$, a
specified state $s_0$ called the \em{start state}, and a specified subset of
states called the \em{accept states}, such that each state has exactly one
outgoing edge labelled with each letter of $\cal A$. A \em{failure state} is any
state which is not an accept state. A \em{path} in $M$ is always a directed
path. Concatenation of paths is denoted by juxtaposition. If $\ell = \ell_\pi$ is
the length of $\pi$, then the states of $M$ visited by $\pi$ are denoted
$\pi[0],\ldots,\pi[\ell]$, and the subpath from $\pi[s]$ to $\pi[t]$ is denoted
$\pi[s,t]$. Reading off the letters on the edges of $\pi$ in succession yields a
word $w_\pi = (a_1 \cdots a_\ell)$ where $a_i$ is the label on the edge
$\pi[i-1,i]$. For any word $w$ and any state $s$, there is a unique directed path
$\pi$ starting at $s$ such that $w = w_\pi$; if $s=s_0$ then we denote this path
by $\pi_w$. When circumstances require, we shall also denote $w_\pi$ by $w(\pi)$
and $\pi_w$ by $\pi(w)$. The set of all words $w$ such that $\pi_w$ ends at an
accept state forms a language over $\cal A$ denoted $L(M)$. 
 
A language $L$ over $\cal A$ is called \em{regular} if there exists a finite
state automaton $M$ over $\cal A$ such that $L=L(M)$. We say that $M$ is a
\em{word acceptor} for $L$. 
 
Given a finite state automaton $M$, an \em{accepted path} is any path from the
start state to an accept state. A \em{live state} is any state lying on an
accepted path. A \em{dead end state} is any state such that all arrows pointing
out of that state point directly back into it; note that a dead end state may be a
live state. Any path which begins and ends at live states is called a \em{live
path}; note that all interior states of a live path are live states. A \em{loop}
is a path which begins and ends at the same vertex, so a live loop is a loop
passing over live states only. Given a loop $\pi$, if we write $\pi = \pi_1
\pi_2$ then $\pi_2 \pi_1$ is also a loop, called a \em{cyclic permutation} of
$\pi$.
 
The basic definition of biautomatic groups involves 2-variable languages (see
e.g. \cite{E$\ldots$, p. 24} or \cite{GS, p. 135}). For our purposes the
equivalent geometric definition of biautomatic groups will suffice,
\cite{E$\ldots$, lemma 2.5.5}, so we shall not use 2-variable languages.
 
Consider a group $G$, an alphabet $\cal A$, and a map $\cal A \to G$. This
induces a monoid homomorphism $\cal A^* \to G$ denoted $w \to \overline w$.
Given $a \in \cal A$, we often omit the overline and consider $a$ as an element
of $G$, even if $\cal A \to G$ is not injective. Thus, abusing terminology,
$\cal A$ is called a \em{generating set} for $G$ if $\cal A^* \to G$ is onto.
Also, any language $L \subset \cal A^*$ that maps onto $G$ is called a set of
\em{normal forms} for $G$. If $M$ is a finite deterministic automaton over $\cal
A$, then for any path $\pi$ in $M$ the group element $\overline w_\pi$ is also
denoted $\overline \pi$.
 
Given a generating set $\cal A$ for $G$, for each $g \in G$ we define the
\em{word length} of $g$ to be $|g| = \Min\{\ell_w \suchthat \overline w = g\}$,
and we define the \em{word metric} on $G$ by $d(g,h) = |g^\inverse h|$.
 
A \em{biautomatic structure} for $G$ consists of a generating set $\cal A$ for
$G$, and a set of normal forms $L \subset \cal A^*$ for $G$, with the following
properties:
\roster
\item $L$ is a regular language
\item There exists a constant $K \ge 0$ such that for each $v,w \in L$ and each
$a \in \cal A \union \{ \epsilon \}$, if $\overline v = \overline w a$ then for
all $t \ge 0$, 
$$ d(\overline v(t), \overline w(t)) \le K
$$
and if $a\overline v = \overline w$ then for all $t \ge 0$,
$$ d(a \overline v(t),\overline w(t)) \le K
$$
\endroster
The constant $K$ is called a \em{two-way fellow traveller constant} for the
biautomatic structure $L$ (to contrast with an automatic structure, in which only
the first inequality is required). As a consequence, for each $v,w \in L$ and any
words $\mu,\nu \in \cal A^*$, if $\overline \mu \, \overline v = \overline w \,
\overline \nu$, then
$$d(\overline \mu \, \overline v(t), \overline w(t)) \le K(|\mu| + |\nu|)
$$
for all $t \ge 0$.
 
We shall need the result of \cite{E$\ldots$, theorem 2.5.1} that any biautomatic
structure on a group $G$ has a sublanguage which is a biautomatic structure with
uniqueness, meaning that each element of $G$ has a unique normal form.
 
Now we review several results of \cite{GS} concerning subgroups of biautomatic
groups; these results will be used without comment in what follows. Fix a
biautomatic structure $L$ on $G$. A subgroup $H \subgroup G$ is called
\em{rational} if the language $\{ w \in L \suchthat \overline w \in H \}$ is
regular. If this is so, then $H$ is a biautomatic group \cite{GS, theorem 3.1 and
theorem 2.2}. The centralizer of a subset $S \subset G$ is denoted $\cal C_S$;
and $\cal C_G$, the center of $G$, is specially denoted $\cal ZG$. If $G$ is
biautomatic and $S$ is a finite set or a finitely generated subgroup, then $\cal
C_S$ is rational \cite{GS, proposition 4.3}; thus, the subgroups $\cal C_S$,
$\cal ZG$ and $\cal Z\cal C_S$ are rational, and it follows that all these
subgroups are biautomatic. 
 
Fix a short exact sequence $$1 \to Z \to G \to H \to 1$$ were $Z = \left< z
\right>$ is an infinite cyclic, central subgroup of $G$. Fix a generating set
$\cal A$ for $G$. Note that $\cal A$ projects to a generating set for $H$ as
well, and the projection map $G \to H$ does not increase the word metric. Let $L$
be a biautomatic structure with uniqueness for $G$ over $\cal A$. Let $M$ be the
word acceptor automaton for $L$. 
 
\subhead A biautomatic structure for $G/Z = H$ \endsubhead
 
Define a \em{central loop} in the automaton $M$ to be any live loop representing
an element of the center $\cal ZG$. We consider two central loops to be the same
if they are cyclic permutations of each other. 
 
\proclaim{Simplicity Lemma (cf. \cite{NS, lemma 3.1})} Let $\gamma$ be a central
loop in $M$. Then $\gamma$ is an iterate of a simple loop in $M$, and every other
simple loop in $M$ is disjoint from $\gamma$. \endproclaim
 
\demo{Proof} If $\gamma$ is not an iterate of a simple loop, then after
cyclically permuting $\gamma$, there exist loops $\mu, \nu$ with the same
initial state as $\gamma$, such that $\gamma = \mu \nu \ne \nu \mu$.
Since $\gamma$ is central, it follows that $\overline \nu \, \overline \mu =
\overline \mu^\inverse \overline \gamma \, \overline \mu = \overline \gamma =
\overline \mu \, \overline \nu$. Choose an accepted path $\pi = \pi_1 \pi_2$
concatenated at the common initial state of the loops $\mu, \nu, \gamma$. Then
$\pi_1 \mu \nu \pi_2$ and $\pi_1 \nu \mu \pi_2$ are distinct accepted paths
representing $\overline \pi \, \overline \gamma$, violating uniqueness of $L$.
 
If there is another simple loop $\gamma'$ in $M$ intersecting $\gamma$, then after
cyclic permutations we may assume that $\gamma$ and $\gamma'$ have the same
base vertex. Since $\gamma$ is central then $\overline \gamma \, \overline \gamma'
= \overline \gamma' \, \overline \gamma$, but $\gamma \gamma' \ne \gamma' \gamma$.
Now proceed as above. \qed\enddemo
 
A central loop is \em{primitive} if it is not an iterate of a shorter central
loop. Note that a primitive central loop does not have to be a simple
loop in $M$. A path $\pi$ in $M$ is said to be \em{compatible with} a set of
primitive central loops $\{ \gamma_1,\ldots,\gamma_I \}$ if $\pi$ intersects each
$\gamma_i$. The set $\{ \gamma_1,\ldots,\gamma_I \}$ is \em{live} if it is
compatible with an accepted path. Define a \em{central cycle} in $M$ to be any
formal linear combination with positive integer coefficients of a live set of
central loops, $c = n_1 \gamma_1 + \cdots + n_I \gamma_I$. The element $\overline
c$ is defined to be $\overline \gamma_1^{n_1} \cdots \overline \gamma_I^{n_I}$.
We say that $c$ is \em{primitive} if $c$ is not a multiple of any other central
cycle, i.e. if the integers $\{ n_1,\ldots,n_I\}$ are relatively prime. 
 
If a path $\pi$ is compatible with a central cycle $c = n_1 \gamma_1 + \cdots +
n_I \gamma_I$, then we may combine $\pi$ and $c$ into a well-defined path as
follows. Choose $t_1,\ldots,t_I$ so that $t_i$ is the minimal integer with
$\pi(t_i) \in \gamma_i$. Since these numbers are distinct by the \tit{Simplicity
Lemma}, we may reindex so that $t_1 < t_2 < \cdots < t_I$. Now take a cyclic
permutation of $\gamma_i$ so that it is based at the point $\pi[t_i]$; this gives
a well-defined loop, since $\gamma_i$ is an iterate of a simple loop. Then the
path 
$$\pi \!*\!* c = \pi[0,t_1] * \gamma_1^{n_1} * \pi[t_1,t_2] * \cdots *
\pi[t_{I-1},t_I] * \gamma_I^{n_I} * \pi[t_I,\ell_\pi]
$$ 
is well-defined. If $\pi$ is an accepted path then $\pi \!*\!* c$ is an accepted
path representing $\overline \pi \, \overline c$. We say that a path $q$
\em{contains} the central cycle $c$ if there exists a path $\pi$ compatible with
$c$ such that $q = \pi \!*\!* c$.
 
A subset of an abelian group is \em{linearly independent} if the identity cannot
be expressed as a non-trivial integer linear combination of elements in the
set. Note that a linearly independent set cannot contain torsion elements.
 
\proclaim{Independence Lemma (cf. \cite{NS, p. 451})} If $\{
\gamma_1,\ldots,\gamma_I \}$ is a live set of central loops, then $\{ \overline
\gamma_1, \ldots, \overline \gamma_I \}$ is a linearly independent subset of
$\cal ZG$. \endproclaim
 
\demo{Proof} Let $\pi$ be any accepted path compatible with $\{ \gamma_1, \ldots,
\gamma_I\}$. If the lemma is false, there is an equation with positive integer
exponents of the form
$$\overline \gamma_{i_1}^{m_{i_1}} \cdots \overline \gamma_{i_A}^{m_{i_A}}
= \overline \gamma_{j_1}^{n_{j_1}} \cdots \overline \gamma_{j_B}^{n_{j_B}}
$$
where $\gamma_{i_a} \ne \gamma_{j_b}$ for $1 \le a \le A$, $1 \le b \le B$. Let
$c_1, c_2$ be the central cycles given by the two sides of this equation, e.g.
$c_1 = m_{i_1} \gamma_{i_1} + \cdots + m_{i_A} \gamma_{i_A}$. Then $\pi \!*\!*
c_1$ and $\pi \!*\!* c_2$ are distinct accepted paths representing the same
element of $G$, contradicting the uniqueness property for $L$.
\qed\enddemo
 
Now define a $Z$-cycle to be a central cycle representing an element of $Z$. If
the element represented is $z^n$ with $n \ge 1$ then the $Z$-cycle is said to be
\em{positive}.
 
\proclaim{Uniqueness Corollary} There are only finitely many primitive $Z$-cycles,
and an accepted path can be compatible with at most one of them. \endproclaim
 
\demo{Proof} There are only finitely many live sets of central loops, and by the
\tit{Independence lemma} each one has at most one primitive linear combination
representing an element of $Z$. If an accepted path $p$ is compatible with two
distinct primitive $Z$-cycles, then those two cycles taken together give a live
set of central loops which forms a linearly dependent subset of $\cal ZG$,
contradicting the \tit{Independence lemma}. \qed\enddemo
 
Define a sublanguage $L_H \subset L$ to consist of all words $w \in L$ such that
$w$ is compatible with some positive $Z$-cycle but $w$ contains no $Z$-cycle. 
 
We shall prove that $L_H$ projects to a biautomatic structure on $H$. The proof
proceeds in three steps: step 1 is to prove that $L_H$ is a regular language;
step 2 is to prove that each coset of $Z$ is represented by some element of
$L_H$; and step 3 is the two-way fellow traveller property.
 
\subhead Step 1: regularity of $L_H$ \endsubhead
 
In one special case the proof of regularity is particularly simple. Namely,
suppose that each $Z$-cycle is actually a simple loop of length 1 in the automaton
$M$. Now delete each such loop, replacing it by an edge leading to a dead end
failure state. The resulting automaton is a word acceptor for $L_H$. 
 
In general, a $Z$-cycle may be a linear combination of non-simple loops. We show
that $L_H$ is regular by reformulating the definition of $L_H$ as a regular
predicate, and then applying \cite{E$\ldots$ theorem 1.4.6}. For any central loop
$\gamma$ in $M$, the language $L_\gamma = \{ w \suchthat \pi_w \intersect \gamma
\ne \emptyset \}$ is regular, because we may alter $M$ by turning each state
lying on $\gamma$ into a dead end accept state, and the new automaton
recognizes $L_\gamma$. Also, the set of words $L^+_\gamma = \{ w \suchthat \pi_w$
contains $\gamma \}$ is regular, because we may modify $M$ by keeping track not
only of the state in $M$ visited by $\pi_w(t)$, but also of the longest subpath of
a cyclic permutation of $\gamma$ traversed by $\pi_w(t)$; this can clearly be done
with a finite state automaton. 
 
For each primitive $Z$-cycle $c = n_1 \gamma_1 + \cdots + n_P \gamma_P$, the
language $L_c = \{ w \suchthat \pi_w$ is compatible with $c \}$ is the same as
$L_{\gamma_1} \intersect \cdots \intersect L_{\gamma_P}$, hence is regular.
Similarly, the language $L^+_c = \{w \suchthat \pi_w$ contains $c \}$ is the same
as $L^+_{\gamma_1^{n_1}} \intersect \cdots \intersect L^+_{\gamma_P^{n_P}}$ hence
is regular.
 
Finally, let $c_1,\ldots,c_N$ be the finite list of all primitive, positive
$Z$-cycles. By the \tit{Uniqueness corollary}, if $w \in L_{c_1}$ then the only
possible $Z$-cycle that $\pi_w$ may contain is $c_1$. Thus, 
$$ L_H = L \intersect \left[ \left( L_{c_1} \intersect \lnot L^+_{c_1} \right)
\union \cdots \union \left( L_{c_N} \intersect \lnot L^+_{c_N} \right) \right] $$
so $L_H$ is regular.
 
\subhead Step 2: $L_H$ represents each coset of $Z$ in $G$ \endsubhead
 
For this argument, fix an element $g \in G$. We must show that the coset $gZ$ is
represented by some word in $L_H$. The proof will depend on the properties of
the Neumann-Shapiro triangulation of the boundary of an automatic structure on
an abelian group.
 
First we make a reduction: it suffices to construct a word $w \in L$ representing
$gZ$ such that $\pi_w$ contains some positive $Z$-cycle $c$. For then we may
write $\pi_w = \pi_v \!*\!* \, c^m$ where $m$ is as large as possible. Then $\pi_v$
does not contain the $Z$-cycle $c$, and yet $\pi_v$ is compatible with $c$, so by
the \tit{Uniqueness corollary} $\pi_v$ does not contain any other $Z$-cycle.
Hence $v \in L_H$ and $v$ represents $gZ$.
 
Let $C = \cal Z \cal C_g$, and note that $C$ contains both $\cal ZG$ and $gZ$. We
review the biautomatic structure on $C$ induced by that on $G$. Let $L_C \subset
L$ be the regular sublanguage of words $w$ with $\overline w \in C$. From the
proof of \cite{GS, theorem 3.1} it follows that there is a generating set $\cal
B$ for $C$, a biautomatic structure $L'$ for $C$ over $\cal B$, and a map $\cal B
\to \cal A$, such that the induced map $\cal B^* \to \cal A^*$ restricts to a
surjection from $L'$ to $L_C$. By \cite{E$\ldots$ theorem 2.5.1} we may replace
$L'$ by a sublanguage which is a biautomatic structure with uniqueness for $C$,
hence the map $L' \to L_C$ is a bijection. Let $M'$ be a word acceptor over
$L'$ over $\cal B$. Then we may speak about $Z$-cycles in $M'$.
 
Now we make another reduction. We shall prove that the coset $gZ$ is represented
by an accepted path $\pi$ in $M'$, so that $\pi$ is compatible with some
$Z$-cycle $c' = n_1 \gamma_1 + \cdots + n_I \gamma_I$ in $M'$. 
 
Accepting this for the moment, we use it to complete step 2. Write $\pi = \pi_0
\pi_1 \cdots \pi_I$ so that for each $k \ge 0$, the path $\pi \!*\!* c'{}^k =
\pi_0 \gamma_1{}^{kn_1} \pi_1 \cdots \pi_{I-1} \gamma_I{}^{kn_I} \pi_I$ is
an accepted path representing $gZ$. Under the mapping $\cal B^* \to \cal A^*$, the
word $w(\pi \!*\!* c'{}^k) \in L'$ goes to a word in the language $L$
corresponding to an accepted path $\pi_k$ in $M$, and each $\pi_k$ represents
$gZ$. For each $k$ and for $1 \le i \le I$, let $t^k_i$ be the moment of time at
which $\pi \!*\!* c'{}^k$ completes the loop $\gamma_i{}^{kn_i}$, and let
$s^k_i$ be the state of $M$ which the path $\pi_k$ visits at time $t^k_i$. Thus,
for each $k$ we obtain an $I$-tuple of states in $M$ denoted $S^k =
(s^k_1,\ldots,s^k_I)$. There must be two distinct values $k_1 < k_2$ with $S^{k_1}
= S^{k_2}$. It follows that $\pi_{k_2} = \pi_{k_1} \!\!*\!* c$ for some $Z$-cycle
$c$ in $M$ representing $\overline c'{}^{k_2 - k_1} \in Z$. Thus, $\pi_{k_2}$ is
an accepted path in $M$ representing $gZ$ and containing a $Z$-cycle, as required
to complete step 2.
 
Now we review the result of Neumann-Shapiro, \cite{NS, theorem 1.1}, which
associates to each automatic structure on the abelian group $C$, a simplicial
decomposition of the boundary. While their result is only stated when $C$ is free
abelian, we note that their construction is valid more generally when $C$ is
abelian.
 
Fix an identification $C = \integers^k \directsum F$ for some finite abelian
group $F$. We shall sometimes confuse an element of $C$ with its projection onto
$\integers^k$. Each non-torsion $c \in C$ determines a ray in $\integers^k$ whose
direction is denoted $[c] \in S^{k-1}$. Neumann and Shapiro associate, to a
biautomatic structure $L'$ on $C$, a rational linear ordered simplicial
subdivision $\Sigma$ of $S^{k-1}$, as follows. Each state of the word acceptor
$M'$ lies on at most one simple loop of $M'$ (see \cite{NS, lemma 3.1}, or the
\tit{Simplicity lemma}). Let $\pi$ be a simple live path in $M'$ initiating at the
start state $s_0$. Let $s_1$ be the first state on $\pi$ which lies on a simple
loop, and let $\gamma_1$ be that loop. Inductively, let $s_i$ be the first state
of $\pi$ after $s_{i-1}$ which lies on a simple loop distinct from
$\gamma_{i-1}$, and let $\gamma_i$ be that loop. This induction ends with $s_l$,
and let $s_{l+1}$ be the final state of $\pi$. Note that $\{ \overline \gamma_1,
\ldots, \overline \gamma_l \}$ is a linearly independent set in $C$ (see
\cite{NS, p. 451}, or the \tit{Independence lemma}). We may now define a rational
linear ordered $(l-1)$-simplex in $S^{k-1}$, namely $\sigma_\pi = \bigl<
[\overline \gamma_1],\ldots,[\overline \gamma_l] \bigr>$. Neumann and Shapiro
prove that as $\pi$ varies over all simple paths in $M'$, the collection $\Sigma
= \{ \sigma_\pi \}$ is an ordered simplicial subdivision of $S^{k-1}$. 
 
Since a group element determines a ray, we need to know the relation between that
element and the simplex at infinity which the ray hits. Let $\pi$ be as above,
and let $\pi_i = \pi[s_{i-1},s_i]$, so $\pi = \pi_1 \pi_2 \cdots \pi_{l+1}$.
Define $\pi(n_1,\ldots,n_l) = \pi_1^{\vphantom{n_1}} \gamma_1^{n_1}
\pi_2^{\vphantom{n_2}} \cdots \pi_l^{\vphantom{n_l}} \gamma_l^{n_l}
\pi_{l+1}^{\vphantom{n_l}}$. Note that each element of $C$ is uniquely represented
by a path of the form $\pi(n_1,\ldots,n_l)$, for some simple accepted path $\pi$
and some $n_1,\ldots,n_l \ge 0$. Fix the ``visual'' metric on $S^{k-1}$, where the
distance between two rays is equal to the angle they subtend. Although we cannot
guarantee that the ray $[\overline \pi(n_1,\ldots,n_l)]$ hits the simplex
$\sigma_\pi$, the following lemma says that it comes visually close:
 
\proclaim{Visual Lemma} For each $\epsilon > 0$ there exists a ball $B \subset C$
around the origin such that if $\overline \pi(n_1,\ldots,n_l) \not\in B$ then the
visual distance between $[\overline \pi(n_1,\ldots,n_l)]$ and the point
$[\overline{\gamma_1^{n_1}} \cdots \overline{\gamma_l^{n_l}}] \in \sigma_\pi$ is
smaller than $\epsilon$. \endproclaim
 
\demo{Proof} Let $\delta$ be the length of the longest simple path in $M'$. Since 
$$\overline \pi(n_1,\ldots,n_0)^\inverse \bigl( \overline{\gamma_1^{n_1}} \cdots
\overline{\gamma_l^{n_l}} \bigr) = \overline \pi^\inverse
$$
it follows that $d(\overline \pi(n_1,\ldots,n_l),\overline{\gamma_1^{n_1}} \cdots
\overline{\gamma_l^{n_l}}) \le \delta$. Choose the ball $B$ so large that a path
of length $\le \delta$ in $C$ touching $C-B$ has visual diameter less than
$\epsilon$. The lemma easily follows. \qed\enddemo
 
To complete step 2, let $\Star[z]$ be the union of those simplices of $\Sigma$
that contain $[z]$. Noting that $[z] \in \interior(\Star[z])$, choose $\epsilon$
so small that every point of $S^{k-1}$ within visual distance $2 \epsilon$ of
$[z]$ is contained in $\interior(\Star[z])$. Choose a positive integer $m$ so
large that $[gz^m]$ is within visual distance $\epsilon$ of $[z]$, and so that $g
z^m$ lies outside the ball $B$ given by the \tit{Visual lemma}. Now $g z^m$ is
represented by a path in $M'$ of the form $\pi(n_1,\ldots,n_l)$, for some simple
accepted path $\pi$ and some $n_1,\ldots,n_l \ge 0$ as above. By the \tit{Visual
lemma} it follows that $[\overline{\gamma_1^{n_1}} \cdots
\overline{\gamma_l^{n_l}}] \in \interior(\Star[z])$, hence $[z] \in \sigma_\pi$.
Therefore there is a positive $Z$-cycle $c$ obtained as a linear combination of
$\gamma_1, \ldots, \gamma_l$. This shows that $g z^m$ is represented by an
accepted path $\pi(n_1,\ldots,n_l)$ in $M'$ compatible with the positive
$Z$-cycle $c$, finishing the proof that the coset $gZ$ is represented by the
language $L_H$.
 
\subhead Step 3: The two-way fellow traveller property for $L_H$ \endsubhead
 
To prove this, consider $v,w \in L_H$ and any $a,b \in \cal A \union \{ \epsilon
\}$, and assume that $a \overline v$ and $\overline w b$ are congruent modulo
$H$. Then $a \overline v = \overline w b z^\beta$ for some $\beta$. We shall give
a bound $|\beta| \le B$, where the constant $B$ depends only on $Z$ and on the
biautomatic structure on $G$. If $K$ is a two-way fellow traveller constant for
$L$, it follows that $d(a \overline v(t), \overline w(t)) \le K' = (B|z|+2)K$ for
all $t \ge 0$. Thus, $K'$ is a two-way fellow traveller constant for $L_H$ in
$G$, and so also in $H$. Henceforth, we can and shall assume $\beta \ge 0$. 
 
To give the idea of the proof, we first sketch the case where $Z$ is a rational
subgroup of $G$. Then we can write $\pi_w = \pi_1 \pi_2$, concatenated at a
vertex that lies on a primitive $Z$-loop $\gamma$. For simplicity, suppose that
$\gamma$ represents $z$ itself, not a power (at worst, $\gamma$ represents a
bounded power of $z$). Then $\overline w z^\beta$ is represented by the accepted
path $\pi' = \pi_1 \gamma^\beta \pi_2$. Let $w' = w_{\pi'} = w_{\pi_1}
w_\gamma^\beta w_{\pi_2}$. Let $k$ be the length of $\gamma$, which is bounded
independent of $\gamma$. Since $a \overline v = \overline w' b$ and $v,w' \in L$,
then  $av$ and $w'b$ are fellow travellers. Assuming by contradiction that
$\beta$ is very large, it follows that $v$ has a long subword $v'$ that fellow
travels the subword $w_\gamma^\beta$ of $w'$. Travelling along $v'$, at every
$\nth{k}$ vertex we keep track of two pieces of data: the state of $M$ visited by
$v'$, and the word difference between $v'$ and $w(\gamma)^\beta$. This data takes
values in a finite set, so if $\beta$ is large enough the data is repeated at two
different spots on $v'$. The subword between these two spots traces out a loop in
$M$, because the states are repeated; and this subword represents an element of
$Z$, because the word difference with powers of $w(\gamma)$ is repeated. Thus, we
have shown that $v$ contains a $Z$-loop, contradicting the fact that $v \in L_H$.
This contradiction shows that $\beta$ cannot be too large, completing the sketch
in the case that $Z$ is rational.
 
In the general case, the path $\pi_w$ can be written in the form $\pi_1 \pi_2
\ldots \pi_p \pi_{p+1}$, with $\pi_j$ and $\pi_{j+1}$ concatenated at a vertex
$V_j$, so that there is a primitive central loop $\gamma_j$ based at $V_j$, and
there is a primitive $Z$-cycle $c = n_1 \gamma_1 + \cdots + n_p \gamma_p$
representing $z^\alpha$, where $0 < \alpha \le A$ for some constant $A$ depending
only on the biautomatic structure on $G$. Now write $\beta = q\alpha + r$ for
some integers $q,r \ge 0$ with $r < \alpha$, so $r < A$. Then there is a word $w'
\in L$ such that 
$$\pi' = \pi_{w'} = \pi_w \!*\!* (q \cdot c) = \pi_1^{\vphantom{qn_1}}
\gamma_1^{qn_1} \pi_2^{\vphantom{qn_1}} \cdots \pi_p^{\vphantom{qn_1}}
\gamma_p^{qn_p} \pi_{p+1}^{\vphantom{qn_1}} 
$$ 
is an accepted path representing $\overline w' = \overline w z^{q\alpha}$. It
follows that $a \overline v = \overline w' b z^r$, so $d(a \overline v, \overline
w' b) \le r|z| < A|z|$. Thus, the words $av$ and $w'$ are fellow travellers with
a constant independent of all choices:  
$$d(a \overline v(t), \overline w'(t)) \le K_1 = (A|z| + 2)K
$$
Let $U$ be the ball of radius $K_1$ around the origin of $G$. 
 
Since the automaton $M$ has only finitely many primitive central loops, for each
such loop $\gamma$ there are only finitely many primitive central
loops having a power representing some power of $\overline \gamma$; let
$G_\gamma$ be this set of loops. There is a positive integer $m_\gamma$
such that each loop in $G_\gamma$ has a power representing
$\overline{\gamma^{m_\gamma}}$. 
 
Recalling the primitive $Z$-cycle $c = n_1 \gamma_1 + \cdots + n_p \gamma_p$,
choose the least positive integral multiple $\rho_j$ of each $m_{\gamma_j}$ so
that $\rho_1 \gamma_1 + \cdots + \rho_p \gamma_p$ is a $Z$-cycle. Note that
$\rho_j$ depends only on the primitive $Z$-cycle $c$ and on $j$. In particular,
there is a global bound $\rho_j \le R$ independent of $c$ and $j$.
 
Fix $j=1,\ldots,p$ for the moment. We show that if $\beta$ is sufficiently large,
then $v$ has an infix subword traversing a loop of $M$ that represents
$\overline{\gamma_j^{\rho_j}}$. Let $L_j$ be the length of $\gamma_j^{\rho_j}$.
Factor $\pi'$ as $\pi'_1 \gamma_j^{q n_j} \pi'_2$, and let the corresponding
factorization of $w'$ be $w'_1 w^+ w'_2$ with $w^+ = w(\gamma_j^{q n_j})$. We may
factor $v = v'_1 v^+ v'_2$, so that for $0 \le t \le qn_j/\rho_j$,   $$d( a
\overline v'_1 \overline v^+(t L_j), \overline w'_1 \overline w^+(t L_j)) \le K_1
$$
Let $d_t$ be this word difference, so $d_t \in U$. Let $s_t$ be the state of $M$
at which the word $v'_1 v^+(t L_j)$ terminates.
 
Noting that $q \ge (\beta - A) / A$, then if
$$\beta \ge {A \rho_j (|U| \cdot |M| + 1) \over n_j} + A
$$
it follows that
$$q \ge {\rho_j (|U| \cdot |M| + 1) \over n_j}
$$
so
$$\left\lfloor {q n_j \over \rho_j} \right\rfloor \ge |U| \cdot |M|
$$
In this case there are integers $0 \le t_1 < t_2 \le qn_j/\rho_j$ so
that $d_{t_1} = d_{t_2}$ and $s_{t_1} = s_{t_2}$. It follows that 
$$\overline v^+(t_1)^\inverse \overline v^+(t_2) = \overline{\gamma_j^{(t_2 -
t_1) \rho_j}}
$$
and that this element is represented by a loop contained in $\pi_v$. This loop
must be an iterate of some simple loop $\gamma' \in G_{\gamma_j}$, and
there must be a lower iterate of $\gamma'$ representing
$\overline{\gamma_j^{\rho_j}}$, since $m_{\gamma_j}$ divides $\rho_j$. Hence,
$\pi_v$ contains a loop reprenting $\overline{\gamma_j^{\rho_j}}$. 
 
Therefore, if $\beta \ge AR(|U| \cdot |M| + 1) + A$ then $\pi_v$ contains a
$Z$-cycle representing $\overline{\gamma_1^{\rho_1}} + \cdots +
\overline{\gamma_p^{\rho_p}} \in Z$, contradicting the fact that $v \in L_H$. This
finishes the proof that $L_H$ is a biautomatic structure for $H$.

\enddocument